\documentclass[modern]{aastex702}

\usepackage{amsmath}

\shorttitle{Toeplitz Determinants and Correlation Intervals}
\shortauthors{Erben}

\begin{document}

\title{Toeplitz Determinants and Admissible Correlation Intervals}

\correspondingauthor{Thomas Erben}

\author{Thomas Erben}
\affiliation{Argelander Institut f\"ur Astronomie, Auf dem H\"ugel
  71, D-53121 Bonn}
\email{terben@astro.uni-bonn.de}

\begin{abstract}
For a homogeneous one-dimensional random field, positive
semidefiniteness of finite Toeplitz correlation matrices imposes
non-trivial constraints on admissible correlation coefficients. The
widths of the corresponding admissible intervals are closely related
to determinants of principal Toeplitz submatrices. Using the classical
Desnanot--Jacobi determinant identity, I derive a simple determinantal
representation for the widths of admissible correlation intervals.

As an immediate consequence, I recover the product expressions for
admissible interval widths previously stated by
\citet{SchneiderHartlap2009}. The argument places these relations into
the general framework of classical Toeplitz determinant theory.
\end{abstract}


\section{Introduction}

For a homogeneous one-dimensional random field with normalized
correlation function \(r(x)\), fix a separation \(x>0\) and define

\[
  r_n := r(nx),\quad n\geq 0,
\]
with \(r_0=1\).

The coefficients \(r_n\) cannot be chosen arbitrarily; rather, positive
semidefiniteness of the finite Toeplitz correlation matrices

\[
  A_N=(r_{|i-j|})_{0\le i,j\le N-1}
\]
imposes non-trivial constraints on the admissible values of successive
correlation coefficients.

For fixed \(r_1,\ldots,r_{n-1}\), the admissible interval boundaries
\(r_{n, \mathrm l}\) and \(r_{n, \mathrm u}\) depend on these previously specified
coefficients and \(r_n\) is restricted to

\[
  r_{n, \mathrm l}\le r_n\le r_{n, \mathrm u},
\]
whose width is

\[
  \Delta_n=r_{n, \mathrm u}-r_{n, \mathrm l}.
\]

Determining these admissible intervals amounts to determining when
finite Toeplitz correlation matrices remain positive semidefinite. In
a study of correlation functions of homogeneous random fields,
\citet{SchneiderHartlap2009}, hereafter SH, derived recursive
constraints on the coefficients \(r_n\) and obtained the product
expressions

\begin{equation}
  \Delta_{2m+1} = 2 \left(\prod_{k=1}^{m}F_{2k}\right)
                   \left(\prod_{k=1}^{m}F_{2k-1}\right)^{-1},
\end{equation}
and

\begin{equation}
  \Delta_{2m} = 2 \left(\prod_{k=1}^{m}F_{2k-1}\right)
                 \left(\prod_{k=1}^{m-1}F_{2k}\right)^{-1},
\end{equation}
where

\[
  F_n=(r_{n, \mathrm u}-r_n)(r_n-r_{n, \mathrm l}).
\]

These product formulae are Equations (24) and (25) of SH. SH stated
these expressions and noted that they had been verified symbolically
with Mathematica for \(n\leq 16\), but no general derivation was
given. They subsequently used these expressions to obtain

\begin{equation}
  \Delta_n=2\,\frac{D_n}{D_{n-1}},\qquad D_n=\det(A_n),
\end{equation}
which appears as Equation (38) of their paper.

The purpose of this note is to show that the interval-width formula
above is a direct consequence of the classical Desnanot--Jacobi
determinant identity \citep{Krattenthaler1999}. The product formulae
of SH then follow immediately. In this way, the admissible correlation
intervals are placed into the framework of classical Toeplitz
determinant theory.

\section{The Desnanot--Jacobi identity}

For a square matrix \(M\), the Desnanot--Jacobi identity states
\citep{Krattenthaler1999}

\[
  \det(M)\det(M_{ij}^{ij}) = \det(M_i^i)\det(M_j^j) - \det(M_i^j)\det(M_j^i),
\]
where \(M_{k}^{i}\) denotes the matrix obtained by deleting row $i$ and
column \(k\), and \(M_{kl}^{ij}\) denotes the matrix obtained by deleting
rows \(i,j\) and columns \(k,l\).

Applying this identity to the corner rows and columns of the
correlation matrix \(A_{n+1}\) gives

\begin{equation}
  D_{n+1}D_{n-1} = D_n^2-E_n^2,
  \label{eq:DJ}
\end{equation}
where \(E_n\) denotes the determinant of the off-diagonal minor
obtained from \(A_{n+1}\) by deleting the first row and the last
column. The second off-diagonal minor is the transpose of the first
and therefore has the same determinant. The off-diagonal minor
defining \(E_n\) contains the coefficient \(r_n\) only once, namely in
its lower-left corner. Hence

\begin{equation}
  E_n = a + b\,r_n,
\end{equation}
with constants \(a\) and \(b\) independent of \(r_n\).

Expanding the off-diagonal minor defining \(E_n\) with respect to this
corner element leaves precisely the central Toeplitz submatrix
$A_{n-1}$. Therefore

\begin{equation}
  b=\pm D_{n-1}.
\end{equation}

\section{Proof of Equation (38) of SH}

I first consider the non-degenerate case in which $A_n$ is positive
definite. In particular, \(D_n>0\) and \(D_{n-1}>0\).

In this case, \eqref{eq:DJ} gives

\[
  D_{n+1} = \frac{D_n^2-E_n^2}{D_{n-1}}.
\]

Since \(E_n=a+b r_n\) and neither \(D_n\) nor \(D_{n-1}\) depend on
\(r_n\), \(D_{n+1}\) is a quadratic polynomial in \(r_n\).  Moreover,
its leading coefficient is \(-b^2/D_{n-1}=-D_{n-1}<0\).

Since \(A_n\) is positive definite, the admissible values of \(r_n\)
are determined by the condition \(D_{n+1}\ge 0\). Hence the interval
endpoints are given by the two zeros of this quadratic polynomial.

The zeros are obtained from

\begin{equation}
  D_n^2-E_n^2=0 \quad\Longleftrightarrow\quad E_n=\pm D_n.
  \label{eq:zeros}
\end{equation}

Since \(E_n\) is linear in \(r_n\), \eqref{eq:zeros} yields the two interval
boundaries \(r_{n, \mathrm l}\) and \(r_{n, \mathrm u}\).

The corresponding values of \(E_n\) differ by \(2D_n\) and with \(E_n = a
+ br_n\) the distance between the two solutions \(r_{n, \mathrm u}\)
and \(r_{n, \mathrm l}\) is

\[
  \Delta_n=r_{n, \mathrm u}-r_{n, \mathrm l}=\frac{2D_n}{|b|}.
\]

Thus Equation (38) of SH follows:

\[
  \Delta_n=2\,\frac{D_n}{D_{n-1}}
\]
as \(|b|=D_{n-1}\) for \(D_{n-1}>0\).

In the same non-degenerate case, a similar argument gives the identity

\begin{equation}
  F_n=(r_{n, \mathrm u}-r_n)(r_n-r_{n, \mathrm l}) = \frac{D_{n+1}}{D_{n-1}}.
  \label{eq:Fn}
\end{equation}

Since \(E_n=a+b r_n\), the expression \(D_n^2-E_n^2\) is a
quadratic polynomial in \(r_n\). Its zeros occur at \(r_{n, \mathrm l}\) and
\(r_{n, \mathrm u}\), corresponding to \(E_n=\pm D_n\), and its leading
coefficient is \(-b^2\). Hence

\[
  D_n^2-E_n^2 = b^2(r_{n, \mathrm u}-r_n)(r_n-r_{n, \mathrm l}).
\]

Using \(b^2=D_{n-1}^2\) and \eqref{eq:DJ}, I obtain

\[
  D_{n+1} = D_{n-1}(r_{n, \mathrm u}-r_n)(r_n-r_{n, \mathrm l})
\]
and consequently \eqref{eq:Fn}.

\section{Degenerate case}

I showed Equation (38) of SH only in the non-degenerate case
\(D_{n-1}>0\).

If \(D_{n-1}=0\), the recursive construction has already reached a
singular Toeplitz matrix. Let \(m\) be the first index for which
\(D_m=0\). Then \(A_m\) is positive semidefinite and singular, whereas
the preceding matrices are positive definite. Hence there exists a
non-zero vector \(v=(v_0,\ldots,v_{m-1})\) such that \(A_m v=0\), and
necessarily \(v_0\ne0\).

To see this, suppose \(v_0=0\). Then the remaining vector
\((v_1,\ldots,v_{m-1})\) would be a null vector of the lower-right
\((m-1)\times(m-1)\) block, which is again \(A_{m-1}\), contradicting
the positive definiteness of \(A_{m-1}\).

Consider any positive semidefinite Toeplitz extension containing the
coefficient \(r_{m+s}\), \(s\ge0\). The principal submatrix with row and
column indices \(0,\ldots,m-1,m+s\) has the block form

\[
  \begin{pmatrix}
    A_m & c^{(s)}\\
    (c^{(s)})^T & 1
  \end{pmatrix},
  \qquad
  c^{(s)}=
  \begin{pmatrix}
  r_{m+s}\\
  r_{m+s-1}\\
  \vdots\\
  r_{s+1}
  \end{pmatrix}.
\]

As a principal submatrix of a positive semidefinite matrix, this block is
positive semidefinite. Evaluating its quadratic form on \((v,t)\in
\mathbb{R}^{m+1}\) gives

\[
  v^T A_m v+2t\,v^T c^{(s)}+t^2\ge0
\]
for all \(t\in\mathbb{R}\). Since \(A_m v=0\), this implies
\(v^T c^{(s)}=0\). Hence

\[
  \sum_{j=0}^{m-1}v_j r_{m+s-j}=0
\]
for \(s\ge0\). Since \(v_0 \ne 0\), this equation determines each new
coefficient \(r_{m+s}\) uniquely from previous ones. Thus all subsequent
admissible intervals collapse to a single point.

The preceding discussion shows that once a singular Toeplitz matrix
occurs, all subsequent admissible intervals have zero width. In this
note, Equation (38) of SH is therefore understood as

\[
  \Delta_n=
  \begin{cases}
    2D_n/D_{n-1}, & D_{n-1}>0,\\
    0, & D_{n-1}=0.
  \end{cases}
\]

\section{Recovery of Equations (24) and (25) of SH}

The initial values are \(D_1=1\) and \(D_2=F_1\). The latter follows
from \(r_{1, \mathrm l}=-1\), \(r_{1, \mathrm u}=1\) and
\(F_1=(1-r_1)(1+r_1)=1-r_1^2=D_2\). Starting from these initial
values, repeated application of \eqref{eq:Fn} yields

\[
  D_{2m+1} = \prod_{k=1}^{m}F_{2k},
\]
and

\[
  D_{2m} = \prod_{k=1}^{m}F_{2k-1}.
\]

For odd \(n=2m+1\), Equation (38) gives

\[
  \Delta_{2m+1} = 2\frac{D_{2m+1}}{D_{2m}} = 2
    \left(\prod_{k=1}^{m}F_{2k}\right)
    \left(\prod_{k=1}^{m}F_{2k-1}\right)^{-1}.
\]

For even \(n=2m\), one obtains

\[
\Delta_{2m} = 2 \frac{D_{2m}}{D_{2m-1}} = 2
  \left(\prod_{k=1}^{m}F_{2k-1}\right)
  \left(\prod_{k=1}^{m-1}F_{2k}\right)^{-1}.
\]

Note that on the degenerate boundary the quotient products are
understood together with the preceding section: once a singular
Toeplitz matrix occurs, all subsequent interval widths vanish.

These are Equations (24) and (25) of SH.

\section{Discussion}

The product formulae of SH follow directly from the Desnanot--Jacobi
identity. The determinant identity provides a direct derivation of the
interval-width formula, from which the product expressions follow
immediately. Although motivated by the correlation-function
constraints studied by SH, the proof itself uses only positive
semidefiniteness of finite Toeplitz matrices and the Desnanot--Jacobi
identity. It therefore applies in the more general setting of positive
semidefinite Toeplitz matrices and places the admissible-interval
relations into the framework of classical Toeplitz determinant theory.

\begin{acknowledgments}

  I thank Matthias Bartelmann and Peter Schneider for helpful comments
  that significantly improved the manuscript.

  I thank Peter Schneider for making me aware of his work while I was
  looking for examples of research projects to which computer algebra
  systems had made a significant contribution. I also thank him for
  clarifications and discussions on his paper.

  I acknowledge the use of OpenAI's ChatGPT as an interactive tool
  during the exploration of the mathematical structure and for
  assistance with drafting and revising parts of the manuscript. I
  take full responsibility for all mathematical arguments and
  conclusions.

\end{acknowledgments}

\bibliographystyle{aasjournalv7.1}
\bibliography{references}

\begin{thebibliography}{}
\expandafter\ifx\csname natexlab\endcsname\relax\def\natexlab#1{#1}\fi
\providecommand{\url}[1]{\href{#1}{#1}}
\providecommand{\dodoi}[1]{doi:~\href{http://doi.org/#1}{\nolinkurl{#1}}}
\providecommand{\doeprint}[1]{\href{http://ascl.net/#1}{\nolinkurl{http://ascl.net/#1}}}
\providecommand{\doarXiv}[1]{\href{https://arxiv.org/abs/#1}{\nolinkurl{https://arxiv.org/abs/#1}}}

\bibitem[{C. Krattenthaler(1999)Krattenthaler}]{Krattenthaler1999}
Krattenthaler, C. 1999, \bibinfo{title}{Advanced determinant calculus,} S\'em.
  Lothar. Combin., 42, Art. B42q, 67 pp. (electronic).
\newblock \doarXiv{math/9902004}

\bibitem[{P. Schneider {\&} J. Hartlap(2009)Schneider \&
  Hartlap}]{SchneiderHartlap2009}
Schneider, P., \& Hartlap, J. 2009, \bibinfo{title}{{Constrained
  correlation functions},} \aap, 504, 705, \dodoi{10.1051/0004-6361/200912424}

\end{thebibliography}

\end{document}